\documentclass{amsart}
\usepackage{amsthm}
\usepackage{amsfonts,amssymb,amsxtra}
\usepackage{xypic}
\usepackage{eufrak}
\usepackage{latexsym}
\usepackage[latin1]{inputenc}
\usepackage{epsf}
\newcommand{\el}{\par \mbox{} \par \vspace{-0.5\baselineskip}}
\newcommand{\goth}[1]{\EuFrak{#1}}
\newcounter{amoi}
\newtheorem{theo}{Th\'eor\`eme}
\newtheorem{leme}{Lemme}
\newtheorem{prop}{Proposition}
\newtheorem{coro}{Corollaire}
\newtheorem{defi}{D\'efinition}

\newcommand{\noi}{\noindent}
\newenvironment{rem}{\noi {\el \noi \bf Remarque}}{ \el }
\newenvironment{dem}{\noi {\el \noi \bf D\'emonstration}}{\hfill $\Box$ \el}

\def\N{{\bf N}}
\def\Z{{\bf Z}}
\def\C{{\bf C}}

\def\today{\number\day \space\ifcase\month\or
 Janvier\or F{\'e}vrier\or Mars\or Avril\or Mai\or Juin\or
 Juillet\or Ao{\^u}t\or Septembre\or Octobre\or Novembre\or D{\'e}cembre\fi
 \space\number\year}

\title[Cônes nilpotents des super algèbres de Lie orthosymplectiques]{Cônes nilpotents des super algèbres de Lie orthosymplectiques}

\author[Caroline Gruson]{Caroline Gruson}
\address{Institut Elie Cartan, UMR 7502 du CNRS, 
Faculte des sciences,
Universit{\'e} Henri Poincar{\'e} (Nancy 1)
BP 239
54506 VANDOEUVRE-les-Nancy Cedex }
\email{gruson@iecn.u-nancy.fr}
\author[S{\'e}verine Leidwanger]{S{\'e}verine Leidwanger}
\address{Institut de math{\'e}matiques de Jussieu, UMR 7586, Equipe de
  Th{\'e}orie 
des Groupes
Case 7012
2 place Jussieu
F-75251 Paris Cedex 05 }
\email{leidwang@math.jussieu.fr}
\date{}
\begin{document}
\keywords{Super algèbres de Lie, Orbites nilpotentes, Désingularisation des cônes nilpotents}
\subjclass{17BXX,14LXX}


\maketitle
\section{Introduction}

Soient $m$, $n$ deux entiers positifs et soit 
$\goth{ g}=\goth{osp}(m,2n)=\goth{ g}_0\oplus \goth{ g}_1$ la super
 algèbre de Lie orthosymplectique correspondante.
On note $ad : \goth{ g}\rightarrow End(\goth{ g})$ l'action adjointe. 
Soit ${\mathcal N_1}$ (resp. ${\mathcal N_0}$) le cône formé des 
$X\in \goth{osp}(m,2n)_1$, (resp. $X\in \goth {osp}(m,2n)_0$) 
tel que $ad(X)$ est un élément nilpotent de 
$End(\goth{ g})$: c'est le nilcône impair (resp. le nilcône pair).

Depuis les travaux de Springer ({\bf [Sp1], [Sp2]}), il est bien connu
que les orbites nilpotentes d'une alg\`ebre de Lie semi-simple $\goth 
a$ sont intimement li\'ees aux repr\'esentations irr\'e\-duc\-tibles
du groupe de Weyl de $\goth a$ et l'un des \'el\'ements de cette
construction de Springer est la d\'esingularisation du c\^one
nilpotent de $\goth a$ par un fibr\'e vectoriel au dessus de la
vari\'et\'e des sous-alg\`ebres de Borel de $\goth a$.

Dans cet article, nous nous intéressons  au cône nilpotent impair ${\mathcal N_1}$ de 
$\goth{ g}$.
D'une part nous étudions les orbites nilpotentes qui le constituent et leurs liens
avec les orbites nilpotentes du cône nilpotent pair ${\mathcal N_0}$, liens obtenus par l'intermédiaire
de l'application $\kappa$, définie par 
$ {\goth{ g}}_1 \rightarrow {\goth{ g}}_0, \  X\mapsto \frac{1}{2}[X,X]$. 
Nous donnons une nouvelle description de la paramétrisation de ces orbites obtenue par Kraft et Procesi dans {\bf[KP]}, 
la nôtre étant fondée sur la paramétrisation des orbites nilpotentes paires (Proposition \ref{pro1}).
Nous étudions ensuite la relation d'ordre (partiel) d'inclusion sur 
les (Zariski) adhérences des orbites nilpotentes impaires en utilisant
 les résultats correspondants de Ohta [{\bf Oh}], Djokovi\'c et Litvinov {\bf[DL]} sur les orbites 
 nilpotentes réelles classiques (Proposition \ref{pro2}).
 
Cette relation d'ordre a pour corollaire l'irréductibilité des fibres de la restriction 
de $\kappa$ à ${\mathcal N_1}$ dont nous donnons une description (Proposition \ref{irr}).
Remarquons que la fibre de $0$ est exactement le cône autocummutant dont 
la géométrie est étroitement reliée à la connaissance de la 
cohomologie de $\goth{ g}$ ({[\bf Gr1]}, {[\bf Gr2]}).

D'autre part, nous donnons, dans les cas $\goth{osp}(2n+1,2n)$ et $\goth{osp}(2n,2n)$, une dé\-singu\-larisation de ce 
nilcône, analogue à celle de Springer pour $\mathcal N_0$, sous la forme d'un fibré 
vectoriel au dessus de la variété des drapeaux de la partie paire (Théorème \ref{th1}).

Dans ce but, nous introduisons la notion de sous-algèbre de Borel mixte (Définition \ref{defmix})
et nous montrons que ces sous-algèbres de Borel rencontrent toutes les orbites nilpotentes. 

Nous remercions chaleureusement Laurent Gruson, Frédéric Han et Nicolas Perrin pour de nombreuses discussions et Michel Duflo pour ses remarques.

\section{Notations, définitions, rappels}

Le corps de base est celui des nombres complexes.

Soient $m\geq 1, n\geq 1$ deux entiers. Soit $V=V_0\oplus V_1$, un espace vectoriel
$\Z/2\Z$-gradué de dimension $m+\varepsilon 2n$ ($dim V_0=m$, $dim V_1=2n$), muni d'une forme 
bilinéaire orthosymplectique non d\'eg\'en\'er\'ee $B$ : on a $B_{|V_0\times V_0}$ est 
symétrique,
$B_{|V_1\times V_1}$ est alternée et $B_{|V_0\times V_1}$ et $B_{|V_1\times V_0}$ sont 
nulles.

On choisit des bases de $V_0$ et de $V_1$ comme suit : si $m=2p$  on prend une base $e'_1,\ldots, e'_{2p}$ 
de $V_0$ avec $B(e'_i,e'_j)=\delta_{i, 2p+1-j}$, si $m=2p+1$ une base $e_0, \ldots e_{2p}$ 
avec $B(e_i,e_j)=\delta_{i, 2p-j}$.
Soit $f_1,\ldots, f_{2n}$ une base de $V_1$  vérifiant 
$B(f_i,f_j)=\delta_{i,2n+1-j}$ si $i\leq j$, $B(f_i,f_j)=-\delta_{i,2n+1-j}$ si $i>j$.
On identifie les endomorphismes de $V$ avec les matrices correspondantes dans ces bases.

On considère la super algèbre de Lie $ \goth{osp}(m,2n)$  constituée des matrices des 
endormorphismes de $V$ qui respectent $B$.
  
On a $\goth{g}=\goth{osp}(m,2n)={\goth g}_0\oplus {\goth g}_1$ avec $\goth{g}_0=\goth{o}(m)\times \goth{sp}(2n)$ et $\goth{g}_1\simeq V_0\otimes V_1$.
Remarquons qu'un élément 
  de $\goth{osp}(m,2n)$ est une matrice par blocs de la forme $\left(\begin{array}{cc}
                         a & u^*\\
                         u & b 
                          \end{array}\right)$, où  $a\in \goth{o}(m)$,
$b\in\goth{sp}(2n)$, $u\in Hom(V_0,   V_1)$ et $u^*$ se déduit de $u$ par 
$$ V_0\xrightarrow{u} V_1 \underset{u^*}{\underrightarrow{\xrightarrow[\sim]{alt} V_1^* \xrightarrow{{}^t u} V_0^* \xrightarrow[\sim]{sym}} V_0}. $$

Soit  $G_0=SO(m)\times SP(2n)$, c'est un groupe algébrique complexe connexe.     
On connaît par {\bf[Vu]} l'anneau $S(\goth{g} _1^*)^{G_0}$, c'est un anneau de polynômes.
On notera ${\mathcal N_1}$ le lieu d'annulation de tous les éléments homogènes non constants de $S(\goth{g} _1^*)^{G_0}$, c'est le cône nilpotent impair (ou nilcône impair de ${\goth g}$).
Remarquons que, selon la terminologie de Mumford, c'est le lieu instable du ${\goth g}_0$-module ${\goth g}_1$ et selon la terminologie de Kac dans {\bf[Ka2]}, le nilcône du ${\goth g}_0$-module ${\goth g}_1$. Il est constitué des éléments  $X$ de ${\goth g}_1$ tels que l'endomorphisme $ad(X)$ de ${\goth g}$ est nilpotent.
On note  ${\mathcal N_0}$ le cône nilpotent de la partie paire.

Soit  $\kappa$ l'application de ${\goth g}_1$ dans ${\goth g}_0$ définie par  $\kappa(X)=\frac{1}{2} [X,X]$.

On a $\kappa: \begin{array}[t]{lll}
 {\mathcal N_1}& \rightarrow &{\mathcal N_0}\\ 
 u& \mapsto & (u^*\circ u,u\circ u^*).
 \end{array}$
 
On remarque que comme le crochet de ${\goth g}$ restreint à la partie impaire est symétrique, $\kappa$ n'est pas identiquement nulle.

\subsection{Diagrammes}

\ \\

Une partition  $\lambda=(\lambda_1,\ldots, \lambda_s)$ est une suite décroissante 
 d'entiers strictement positifs. Les $\lambda_i$ sont appelées 
{\it parts} de $\lambda$. Le nombre de parts, noté $l(\lambda)$, est  la longueur 
 de $\lambda$.  Si $\sum_1^s  \lambda_i=n$, on dit que $\lambda$ est une partition de $n$.
 On appelle sous-partition toute suite décroissante d'entiers formée d'une sous-suite décroissante de $\lambda$ et éventuellement de $0$.

On note $D(\lambda)$ le diagramme de Young de forme $\lambda$.

\begin{figure}[h!]
\leavevmode
 \epsfxsize=7cm 
\centerline{\epsfbox{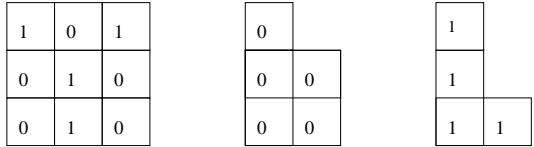}} 
\caption{\label{fig2}$D, d_0, d_1$}
\end{figure}

\begin{defi}\label{def1}  Un diagramme gradué $D$ est un diagramme de Young dans lequel 
chaque case est remplie par un $0$ ou un $1$, de telle manière que les étiquettes alternent sur les  
lignes. Deux diagrammes gradués sont considérés comme égaux si l'on peut passer de l'un à 
l'autre en échangeant des lignes de même longueur.
Un sous-diagramme de $D$ est un diagramme obtenu en effaçant des lignes de $D$.
\end{defi}

\`A un diagramme gradué on associe deux partitions $d_0$ et $d_1$, correspondant 
aux diagrammes obtenus en effaçant les cases $1$ (resp. les cases $0$) dans le
 diagramme gradué et en réordonnant les lignes afin d'obtenir le diagramme 
 d'une partition. On peut voir un exemple sur la figure \ref{fig2}.

Une ligne commençant par une case étiquetée  $0$ (resp. $1$) est dite {\it paire} (resp.
{\it impaire}).

\subsection{Paramétrisation des orbites nilpotentes des algèbres de Lie symplectiques 
et orthogonales}
\ \\

La paramétrisation des orbites nilpotentes  des algèbres de Lie semi-simples complexes 
est décrite dans {[\bf CM]}.

\begin{prop}{\bf [CM]}

\begin{itemize}
\item  Les orbites nilpotentes sous $SO(2n+1)$ de $ \goth{o}(2n+1)$ sont param{\'e}tr{\'e}es par les
partitions de $2n+1$ pour lesquelles les parts paires apparaissent avec
une multiplicit{\'e} paire, on note indifféremment $P_{\goth{o}(2n+1)etiq}$ ou   $P_{\goth{o}(2n+1)}$ l'ensemble de ces partitions. 

On a par exemple  $P_{\goth{o}(5)etiq}=\{(5),(3,1,1),(2,2,1),(1,1,1,1,1)$\}. 
\item Celles de $ \goth{o}(2n)$ sous $O(2n)$ sont paramétrées par les partitions $\lambda$ de $2n$ pour 
lesquelles les parts paires apparaissent avec multiplicité paire, notons $P_{\goth{o}(2n)}$ cet ensemble.
\item Celles de $ \goth{o}(2n)$ sous $SO(2n)$ sont paramétrées par les partitions $\lambda$ de $2n$ pour 
lesquelles les parts paires apparaissent avec multiplicité paire excepté pour les
partitions  formées uniquement de parts paires apparaissant avec des multiplicités
 paires, appelées {\bf très paires}. Dans ce cas on considère les partitions étiquetées $\lambda^I$ et $\lambda^{II}$ correspondant à deux orbites différentes. Les étiquettes $I$ et $II$ étant définies grâce aux diagrammes de Dynkin pondérés {\bf lem 5.3.5 [CM]}. 
 On note $P_{\goth{o}(2n)etiq}$ l'ensemble de ces partitions étiquetées.
 On a par exemple 
 $P_{\goth{o}(8)etiq}=\{(7,1),((5,3),(4^2)^I, (4^2)^{II},(5,1^3),(3^2,1^2),
 (3,2^2,1),(3,1^5),(2^4)^I,$ $ (2^4)^{II},(2^2,1^4), (1^8)\}.$

\item  Celles de $\goth{ sp}(2n)$ sous $SP(2n)$ sont param{\'e}tr{\'e}es par des partitions de $2n$
pour lesquelles les parts impaires apparaissent avec une multiplicit{\'e}
paire, notons $P_{\goth{ sp}(2n)}$ l'ensemble de ces
partitions. On a par exemple $P_{\goth{ sp}(4)}=\{(4),(2,2),(2,1,1),(1,1,1,1)\} $.
\end{itemize}
\end{prop}

On en déduit une paramétrisation des $G_0$-orbites nilpotentes de ${\mathcal N_0}$ de la 
super algèbre de Lie $\goth{ osp}(m,2n)$.
Elles sont paramétrées par des paires de partitions 
$(\lambda,\mu)\in P_{\goth{o}(m)etiq}\times P_{\goth{ sp}(2n)}$.

\subsection{Paramétrisation des $O(m)\times SP(2n)$-orbites nilpotentes de
 $Hom(V_0,V_1)\times Hom(V_1,V_0)$}
\ \\
 
Nous rappelons ici les résultats de Kraft et Procesi ({[\bf KP]}).

\begin{prop}{\bf[KP]}\label{indec}

Les $O(m)\times SP(2n)$-orbites nilpotentes de $\goth{g}_1$ sont paramétrées par des diagrammes gradués 
$D$ formés des diagrammes gradués indécomposables décrits ci-après. On note $O_D$
l'orbite associée à $D$.

Soit $p$ un entier quelconque, les différents diagrammes indécomposables sont
\begin{enumerate}
\item une ligne paire de longueur $4p+1$, $p\geq 0$
\item une ligne impaire de longueur $4p-1$, $p\geq 1$
\item deux lignes paires de longueur $4p-1$, $p\geq 1$
\item deux lignes impaires de longueur $4p+1$, $p\geq 0$
\item deux lignes, l'une paire l'autre impaire de longueur $2p$, $p\geq 1$.
\end{enumerate}
On note $\mathcal D(m,2n)$ l'ensemble des diagrammes gradués formés par 
ces diagrammes gradués indé\-com\-po\-sables.
\end{prop}

\section{L'application $\kappa$}

Dans le cas de la super algèbre de Lie $\goth{ osp}(m,2n)$ avec $m$ 
impair, les $O(m)\times SP(2n)$-orbites nilpotentes de $\goth{g}_1$  et les $G_0$-orbites de $\goth{g}_1$
 sont les mêmes, la proposition précédente nous permet donc de
  paramétrer les orbites nilpotentes de $\goth{ g}_1$.

Dans le cas de $\goth{ osp}(m,2n)$ avec $m$ pair, ce n'est pas le cas. 
Certaines $O(m)\times SP(2n)$-orbites sont l'union de deux $G_0$-orbites
  distinctes dans $\goth{g}_1$ qu'on ne peut différencier par les diagrammes. 

On donne dans un premier temps une nouvelle description des diagrammes apparaissant dans la proposition \ref{indec}.
\begin{prop}\label{pro1} 

 Pour $m\geq 1$, les $O(m)\times SP(2n)$-orbites nilpotentes de $\goth{ osp}(m,2n)$ dans $\goth{g}_1$ 
sont paramétrées par des diagrammes gradués $D$ de forme $\nu$ une partition de $m+2n$
 tels que
\begin{enumerate}
\item le nombre $k$ de lignes de longueur paire de $D$ est pair, et dans $D$ il y a  
$k/2$ lignes paires (resp. impaires) de cette longueur.
\item Soit  $(d_0,d_1)$ le couple de partitions associé à $D$, on a 
$(d_{0}, d_{1})\in P_{\goth{o}(m)}\times P_{\goth{ sp}(2n)}$. 

\end{enumerate}
\end{prop}

\begin{dem} - Soient $D\in \mathcal{ D}(m,2n)$ et ${\bf O}_D$ la $O(m)\times SP(2n)$-orbite associée, la condition $1$ de la proposition 
provient du cas $5$ de la proposition \ref{indec}. 
Soit $X=(u,u^*)\in  {\bf O}_D$. Les diagrammes gradués nous donnent la dimension des noyaux des $X^i$.
Or on a 
$[X,X]=2\left(\emph{} \begin{array}{cc} u^*\circ u & 0 \\ 0& u\circ  u^*\\ \end{array} \right)\in {\mathcal N_0}$. 
Lorsque dans le diagramme gradué on efface les cases étiquetées $1$ (resp. $0$) on ne garde que 
la dimension des noyaux des homomorphismes de $V_0$ dans $V_0$ (resp. $V_1$ dans $V_1$) obtenus 
en restreignant les $X^{2i}$ à $V_0$ (resp. $V_1$). On obtient ainsi une partition indexant 
une (ou deux) orbites de $\goth{o}(m)$ (resp. $\goth{{sp}}(2n)$.
On a donc  $(d_{0}, d_{1})\in P_{\goth{o}(m)}\times P_{\goth{sp}(2n)}$. 

Montrons que 
$(d_0,d_1)\in P_{\goth{o}(m)\times \goth{ sp}(2n)}$ et que deux diagrammes gradués de même 
forme ne peuvent avoir deux couples de partitions associées égaux.

Supposons qu'il existe deux diagrammes $D$ et $D'$ distincts de forme $\nu$ et de couples de 
partitions associées $(\lambda, \mu)$.
Cela signifie que, dans l'un des diagrammes on a une ligne paire de longueur $2k$ (et pas de 
ligne impaire de même longueur), et dans l'autre une ligne impaire de longueur $2k$ (et pas de 
ligne paire de même longueur).
Mais ceci est en contradiction avec le fait que le diagramme gradué est formé d'éléments 
de la proposition \ref{indec}.
 
Il reste à montrer qu'avec un couple de partitions 
$(\lambda, \mu)\in P_{\goth{o}(m)}\times P_{\goth{ sp}(2n)}$ on peut former
un diagramme gradué, et qu'en rajoutant la condition $1$ de la proposition, ce diagramme est 
toujours dans $\mathcal {D}(m,2n)$.

Soit $\lambda\in P_{\goth{o}(m)}$. Elle est formée d'éléments indécomposables de deux sortes
\begin{itemize}
\item cas (0-1) une part de longueur impaire, 
\item cas (0-2) deux parts de longueur paire.
\end{itemize}
Soit $\mu \in P_{\goth{ sp}(2n)}$. Elle est formée d'éléments indécomposables de deux sortes, 
\begin{itemize}
\item cas (1-1) une part de longueur paire,
\item cas (1-2)  deux parts de longueur impaire.
\end{itemize}
\'Etudions tous les couples possibles aboutissant à des diagrammes indécomposables. Pour 
pouvoir construire une ligne de diagramme gradué à partir de ces indécomposables il faut 
que les parts qu'on associe vérifient $\vert \lambda_i-\mu_j\vert\leq 1$.
\begin{itemize}
\item cas (0-1) $\times $cas (1-1): on obtient soit une ligne impaire de longueur $4p-1$, soit une ligne 
paire de longueur $4p+1$ i.e. cas $1$ et $2$ de la proposition \ref{indec}
\item cas (0-2) $\times$ cas (1-2): on obtient soit deux lignes impaires de longueur $4p+1$, soit deux 
lignes paires de longueur $4p-1$, cas $3$ et $4$ de la proposition \ref{indec}.
\item cas (0-1) $\times$ cas (1-2): soit dans $\lambda$ on a deux parts impaires de longueur $2p+1$ et 
on obtient une ligne paire de longueur $4p+2$, une ligne impaire de longueur $4p+2$ i.e 
cas $5$ de la proposition \ref{indec}, soit ce n'est pas le cas et on obtient alors une seule ligne paire
 et ce cas est incompatible avec la condition $1$ de la proposition.
\item cas (0-2) $\times$ cas (1-1): se traite de la même manière que le cas précédent.
\end{itemize}

\end{dem} 

Pour pouvoir ensuite paramétrer les $G_0$-orbites nilpotentes on introduit la définition suivante.
\begin{defi}

Soit $D\in \mathcal{D}(m,2n)$ un diagramme gradué vérifiant les hypothèses de  la proposition précédente. Si $d_0$ est très paire alors à $d_0^I$ (resp. $d_0^{II}$) on associe un diagramme gradué étiqueté $D^I$ (resp. $D^{II}$) sinon on garde $D$. L'ensemble des diagrammes gradués ainsi construit est appelé ensemble des diagrammes gradués étiquetés et noté $\mathcal{D}(m,2n)_{etiq}$.
\end{defi}

\begin{prop}\label{proplab}

Les $G_0$-orbites nilpotentes de $\goth{ osp}(m,2n)$ dans $\goth{g}_1$ sont paramétrées par des diagrammes gradués étiquetés.
\end{prop}

Soit $D=(d_0^i, d_1)$ un diagramme gradué étiqueté, on note $O_D$ la $G_0$-orbite nilpotente telle que $\kappa(O_D)=O_{d_0^i,d_1}$.

\begin{dem} - 
Soit $D\in \mathcal{D}(m,2n)_{etiq}$ et $(d_0^i, d_1)$ son couple de partitions associé.
Si $d_0$ est très paire (i.e. $i=I$ ou $i=II$) alors clairement $O_D$ est unique.
Si $d_0$ n'est pas très paire le résulat découle de la proposition précédente, et du lemme 4.3 de [{\bf KP}].

\end{dem}

Il existe une relation d'ordre (partiel) d'inclusion sur 
les (Zariski) adhérences des $G_0$-orbites nilpotentes impaires.
Soit $D \in {\mathcal D}(m,2n)$, $D'\in {\mathcal D}(m,2n)$, on écrit que  $O_D \leq O_{D'}$ si $O_D
$ est contenue dans l'adhérence de Zariski de $O_{D'}$.
  
Il en existe aussi une sur l'ensemble ${\mathcal D}(m,2n)$ (resp. $\mathcal{D}(m,2n)_{etiq}$).
On peut en trouver une définition dans {\bf[Oh]}( resp. {\bf [DL]}).
\begin{defi}\label{def2}{\bf[Oh]}  \begin{itemize}
\item Soit $D\in {\mathcal D}(m,2n)$, on désigne par  $\underline{D}$ le diagramme gradué obtenu 
à partir de $D$ en effaçant la première colonne.
Soit $k$ un entier supérieur ou égal à $1$.
On définit $D^{(k)}$ par récurrence,  on pose $D^{(0)}=D$
et $D^{(k)}:=\underline{D^{(k-1)}}$.
\item Soient $D,D' \in \mathcal{D}(m,2n)$, on a $D\leq D'$
  si pour $k\geq 1$,
$n_{0}(D^{(k)})\leq n_{0}(D'^{(k)})$ et $n_{1}(D^{(k)})\leq
n_{1}(D'^{(k)})$ 
 où $n_{0}(D)$ (resp . $n_{1}(D)$) est le nombre de 
${0}$ 
(resp. $1$) dans $D$.
\end{itemize}
On note $\Gamma(m, 2n)$ le diagramme de Hasse de cet ordre sur $\mathcal{D}(m,2n)$.
\end{defi}

\`A partir de  $\Gamma(m,2n)$ on  construit un diagramme de Hasse $\Delta(m, 2n)$  sur $\mathcal{D}(m,2n)_{etiq}$
(en suivant {[\bf DL]}) afin de prendre en compte le cas des diagrammes $D$ pour lesquels 
$d_0$ est une partition très paire, et qui indexent donc deux orbites différentes.
Soient $D,D' \in {\mathcal D}(m,2n)$ deux diagrammes différents, on note $D\rightarrow D'$, si $D<D'$ et 
s'il n'existe pas $D''\in {\mathcal D}(m,2n)$ tel que $D< D"< D'$.
Si la $O(m)\times SP(2n)$-orbite paramétrée par $D$ n'est pas connexe on dit que $D$ est {\it instable} 
(sinon il est dit {\it stable}).

\begin{defi}{\bf[DL]}

$\Delta(m, 2n)$ est obtenu à partir de $\Gamma(m,2n)$ en faisant les modification suivantes:
\begin{enumerate}
\item pour chaque paire $(D,D')$ telle que $D\rightarrow D'$ et $D$ ou $D'$ est instable, on efface 
l'arête entre $D$ et $D'$.
\item On remplace chaque noeud $D$ par  deux noeuds, $D^I$ et $D^{II}$.
\item On insère deux arêtes pour chaque arête effacée dans le point 1 de la façon suivante:
 \begin{itemize}
    \item si $D$ est stable et $D'$ instable: on joint $D$ à $D'^{I}$ et à $D'^{II}$,
    \item si $D'$ est stable et $D$ instable: on joint  $D^{I}$ et  $D^{II}$ à $D'$,
    \item si $D$ et $D'$ sont tous deux instables: on joint $D^I$ à $D^{'I}$ et 
    $D^{II}$ à $D^{'II}$.
  \end{itemize}
    \end{enumerate}
\end{defi}

Comme dans le cas des paires symétriques ({\bf[Oh]}, {\bf [DL]}) on a la proposition suivante:

\begin{prop}\label{pro2} Soient $D$ et $D'$ des diagrammes gradués étiquetés. On a $O_D\leq O_{D'}$ si et seulement si  $D\leq D'$ où ($\leq$) a pour diagramme de Hasse $\Delta(m,2n)$.  On a donc  
$\overline{O_{D'}}=\bigcup_{D\leq D'} O_{D}$.
\end{prop}

\begin{dem} -  Le fait que, si $O_D\subset \overline{O}_{D'}$  alors $D\leq D'$, est évident.
La démonstration de l'ordre sur les adhérences des $O(m)\times SP(2n)$-orbites est basée sur les 
résultats de {\bf[Oh]}, celle sur l'ordre des adhérences des $G_0$-orbites sur ceux de
{\bf[DL]}.
On rappelle en annexe les différentes notations liées aux paires symétriques en suivant celles de {\bf [Oh]}.

\begin{leme}  Soient $D$ et $D'$ deux diagrammes gradués tels que 
\begin{itemize}
\item leurs premières colonnes coïncident,
\item $D\leq D'$, où $\leq$ est l'ordre de la définition \ref{def2},
\item $D$ et $D'$ indexent les orbites de type 
(DIII)  ou (CI) (cf annexe). 
\end{itemize}
Soient $O_{(DIII), D}$ et $O_{(DIII), D'}$ (resp. $O_{(CI), D}$ et $O_{(CI), D'}$)\emph{} ces orbites.
Dans ce cas $ D^{(1)}$ et $ D'^{(1)}$ indexent deux $O(m)\times SP(2n)$-orbites $O_{ D^{(1)}}$ et $O_{D^{'(1)}}$
 et on a $ D^{(1)}\leq  D'^{(1)}$.
Si  $O_{(DIII), D}\subset \overline{O}_{(DIII), D'}$ 
(resp. $O_{(CI), D}\subset \overline{O}_{(CI), D'}$) 
on a alors $O_{D^{(1)}}\subset \overline{O}_{ D^{'(1)}}$.  
\end{leme}

\begin{dem} -  On adapte la démonstration du lemme 7 de {\bf [Oh]}.
Les rappels de définition et de notation concernant les paires symétriques sont donnés en annexe.

En effaçant la première colonne de ces éléments indécomposables on retrouve les éléments indécomposables 
indexant les $O(m)\times SP(2n)$-orbites.  
On va démontrer le lemme pour les orbites de type (DIII), l'autre cas se traitant de la même façon.

Soient $W=W_0\oplus W_1$, $U=U_0\oplus U_1$ deux espaces vectoriels munis d'une involution notée 
$s_W$ (resp. $s_U$). 
On munit $W$ d'une forme bilinéaire orthogonale non dégénérée $B_W$ on a : $B_{W_{|W_i\times W_i}}$ est orthogonale, $B_{W_{|W_0\times W_1}}$ et $B_{W_{|W_1\times W_0}}$ sont 
nulles. On munit  $U$ d'une forme bilinéaire 
orthosymplectique non dégénérée notée $B_U$.

\noi On pose : $$L(W,U):=Hom(W,U), 
L^-(W,U):=\{X\in L(W,U)/ s_W X s_U=-X\}.$$

\noi On définit l'adjoint $X^*\in L(U,V)$ de $X$ par $B_U(Xw,u)=B_W(w, X^*u)$ pour $(w,v)\in W\times V$.

\noi On note $K(U)=O(m)\times SP(2n)$, $K(W)=O(W)\cap GL(W_0)\times GL(W_1)$.

\noi On a alors : $K(W)\times K(U)$ agit sur $L^- (W,U)$ par $(g,h)X=gXh^{-1}$.

\noi En appliquant le lemme 11 ({\bf[Oh]}), on peut définir deux morphismes :

 $$ \goth{p}(W)\stackrel{\rho}{\leftarrow}L^-(W,U)
 \stackrel{\pi}{\rightarrow}\goth{ g}_1, 
 \pi(X)=XX^*, \rho(X)=X^*X$$.

\noi On pose $M:=\rho^{-1}(O_{(DIII),D'})$, on a $\rho(\overline{M})=\overline{O}_{(DIII),D'}$ 
({\bf[Oh]}, p206).
De même $\pi(\overline{M})= \overline{O}_{ D'^{(1)}}$.
Puisque $D\leq D'$ on a $O_{(DIII), D}\subset \overline{O}_{(DIII), D'}$. Or on a $\overline{O}_{(DIII), D'}=\rho(\overline{M})$.
Il existe donc $Y\in \overline{M}$ tel que $\rho(Y)\in O_{(DIII), D}$. On a, d'autre part, 
$\pi(\rho^{-1}(O_{(DIII), D}))=O_{ D^{(1)}}$ (lemme 14,(2),{\bf[Oh]}).
On en déduit que $\pi(Y)\in O_{ D^{(1)}}\cap \pi(\overline{M})$ c'est-à-dire $\pi(Y)\in \overline{O}_{ D'^{(1)}}$.
On a donc $O_{D^{(1)}}\subset \overline{O}_{D'^{(1)}}$.
\end{dem}

En jouant sur le fait qu'on peut échanger le rôle des $0$ et des $1$ dans les cas traités par {\bf [Oh]}, 
on remarque que tous les cas  qui restent sont traités dans la table {\sc v} de {\bf [Oh]}.

Il reste à présent à démontrer que si $D\leq D'$, où $\leq$ a pour diagramme de Hasse $\Delta(m,2n)$,
alors 
les $G_0$-orbites associées vérifient $O_D\subset \overline{O}_{D'}$.

Pour cela on utilise les résultats de {\bf [DL]}.
Soient $D$ et $D'$ deux diagrammes gradués instables tels que $D'< D$ dans $\Gamma(m,2n)$. Si, pour tout $D"$ tel que  $D'<D''<D$,
$D''$ est instable, on dit que $(D,D')$ est une {\it paire pure}.

\begin{leme}
Si $(D,D')$ est une paire pure alors $d_0=d_0'$.
\end{leme}

\begin{dem} -
On est dans le cas où $m=2p$.
  
Soit $D$ un diagramme gradué tel que $d_0$ est très paire.

On va démontrer le lemme en montrant qu'il est vrai pour les indécomposables.
On vérifie facilement que les seuls indécomposables possibles avec $d_0$ très paire correspondent aux cas (3),(4),(5) de la proposition \ref{indec}.
On  traite en détail le cas (3), le raisonnement étant similaire pour les autres cas.

\noi On construit tous les diagrammes possibles $D'<D$ tels que $d_0'<d_0$ et $d_0'$ est très paire.

\noi On peut les voir sur la figure \ref{inpair}.

Soit $D"=(d_0",d_1")$ le diagramme formé 
\begin{itemize}
\item d'un indécomposable correspondant au cas (5) obtenu à partir du cas (3) en enlevant deux cases paires, 
\item de deux lignes paires de longueur $1$. 
\end{itemize}
Il est représenté sur la figure \ref{inpair}.

\noi On vérifie aisément que $d_0"$ n'est pas très paire et que $D'<D"<D$.
\begin{figure}[h!]
\leavevmode
 \epsfysize= 8cm 
\centerline{\epsfbox{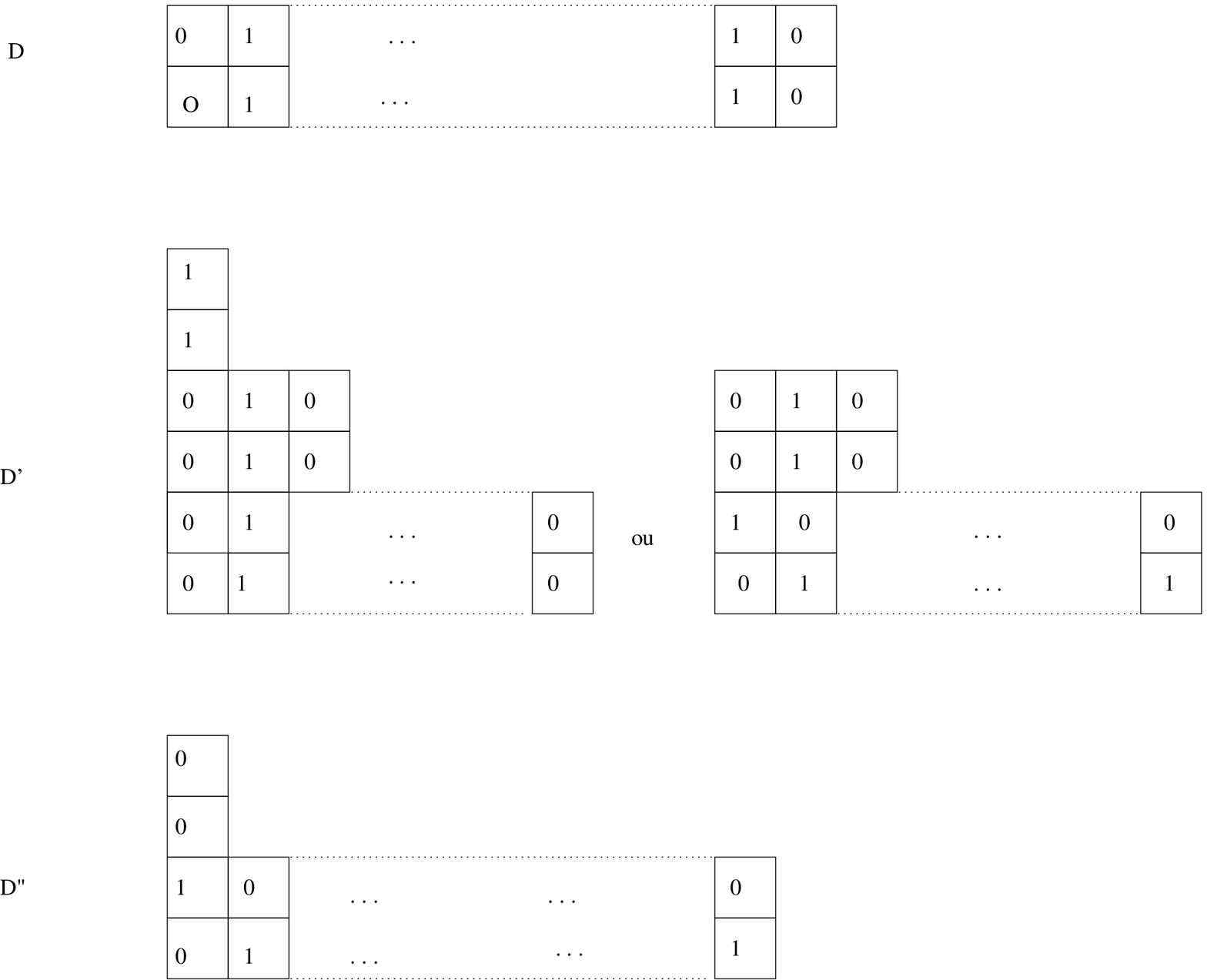}} 
\caption{\label{inpair}}
\end{figure}

On a donc montré que si, $D'<D$ avec $d_0'$ et $d_0$ très paires et $d_0'<d_0$, alors la paire $(D,D')$ n'est pas pure.
\end{dem}
 Il suffit 
de montrer (Théorème 3.7 {\bf[DL]}) que pour une paire pure $(D,D')$, $\overline{O}_{D^I}\cap  O_{D'}=O_{D'^I}$.
Or on a :

$O_{D'^I}\cup O_{D'^{II}}= O_{D'}\subset \overline O_{D}=\overline{O}_{D^I}\cup \overline{O}_{D^{II}}$.

Si $O_{D'^I}\subset \overline{O}_{D^{II}}$  alors, par continuité de $\kappa$, on a $\kappa(O_{D'^I})\subset \kappa(\overline{O}_{D^{II}})$, i.e. $O_{{d'_0}^I, d'_1}\subset \overline{O}_{d_0^{II}, d_1}$, ce qui n'est pas le cas.
On a donc $\overline{O}_{D^I}\cap O_{D'}=O_{D'^I}$. 
\end{dem}

\begin{figure}[h!]
\leavevmode
 \epsfysize= 20cm 
\centerline{\epsfbox{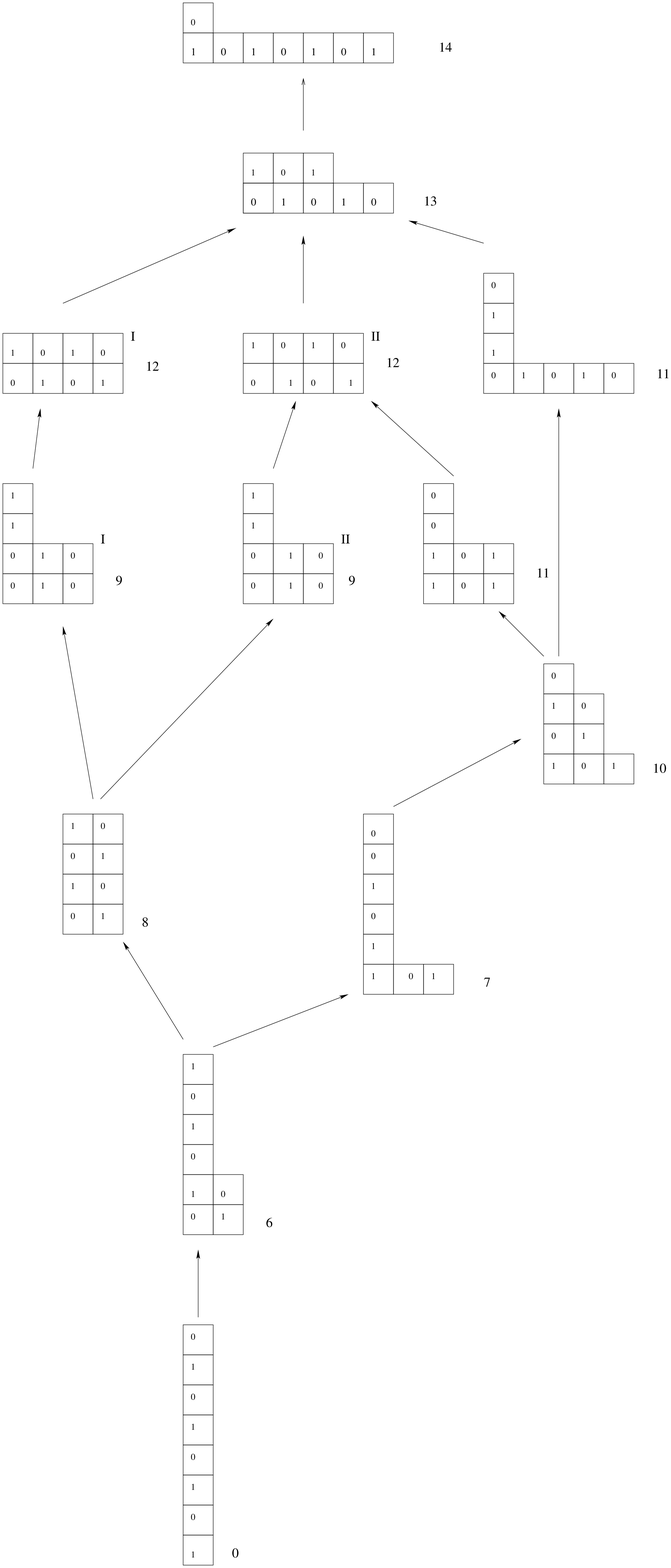}} 
\caption{\label{fig6}}
\end{figure}

On peut voir le diagramme de Hasse de l'ordre sur les $G_0$-orbites nilpotentes impaires
de $\goth{ osp}(4,4)$ sur la  figure \ref{fig6} (les orbites sont représentées par leurs diagrammes gradués étiquetés, on a aussi indiqué la dimension des orbites que l'on calculera dans la suite).

\vskip 1cm

On s'intéresse à présent à l'application $\kappa$. 
On décrit ses fibres au dessus de chaque orbite paire. En utilisant la proposition \ref{pro1} on 
obtient facilement le lemme ci-dessous.

\begin{leme}\label{cor1}  Soit $O_{\lambda,\mu}$ une $G_0$-orbite nilpotente de  $\goth{ g}_0$
\[ \kappa^{-1}(O_{\lambda,\mu})=\left \{ 
              O_{D}\in {\mathcal N_1},  D \in {\mathcal D}(m,2n) \mbox{ tel que
              } 
d_{0}=\lambda, d_{1}=\mu \}  
                     \right. \]

\end{leme}

\begin{rem} -\label{rema}  La fibre est vide la plupart du temps car c'est le cas si dans $\mu$
il n'y a pas de part égale à $\lambda_i-1$, $\lambda_i$ ou $\lambda_i+1$.
\end{rem}
\vfill\eject

\begin{leme}\label{irr} Soient $(\lambda,\mu)\in  P_{\goth{o}(m)etiq}\times P_{\goth{ sp}(2n)}$ telles que $\lambda$ et $\mu$  comportent au plus une sous-partition en commun de la forme $(k,k,k-1,k-1)$, $k\geq1$. Il existe un unique diagramme gradué étiqueté maximal  $D_{\scriptsize{\mbox{ max}(\lambda,\mu)}}$ vérifiant 
$(d_0,d_1)=(\lambda,\mu)$. On a 
   $\kappa^{-1}(\overline{
O}_{\lambda,\mu})=\overline{O}_{D_{\scriptsize{\mbox{ max}(\lambda,\mu)}}}$. En d'autres termes, 
 $\kappa^{-1}(\overline{O}_{\lambda,\mu})$ est irréductible.
\end{leme}

\begin{dem} -  Supposons que $\lambda$ et $\mu$  ne comportent pas de sous-partition en commun de la forme $(k,k,k-1,k-1)$ avec $k>1$. Pour que $\kappa^{-1}(O_{\lambda,\mu})$ contienne plus d'une orbite il faut 
que dans $\lambda$ (resp. $\mu$) il y ait $2k$
(resp. $2l$) 
parts égales à $1$. On remarque alors que si $\lambda$ est  très paire 
$\kappa^{-1}(O_{\lambda^{i},\mu})$ où $i \in\{I,II\}$ 
ne contient qu'une orbite.
On en déduit que si $\kappa^{-1}(O_{\lambda,\mu})$ contient plusieurs orbites alors 
l'ordre sur ces orbites a pour diagramme de Hasse $\Gamma(m,2n)$.
 Il suffit ensuite de montrer que  $D_{\scriptsize{\mbox{ max}(\lambda,\mu)}}$ existe.
 Ceci résulte du fait que l'ensemble
${\bf D_{\lambda,\mu}}=\{ D \in {\mathcal D}(m,2n)/
  (d_0,d_1)=(\lambda,\mu)\}$ est totalement ordonné.

Soit $D_{res}$ (resp. $D_{ind}$) le diagramme obtenu à partir de 
$D\in {\bf D_{\lambda,\mu}}$ en
effaçant toutes les lignes de longueur $l>2$ (resp. de longueur $l\leq 2$).
Tous les  $D_{ind}$ sont les mêmes pour $D\in {\bf  D_{\lambda,\mu}}$. 

Les diagrammes $D_{res}$  possibles sont dessinés sur la figure \ref{deca}.
\begin{figure}[h!]
\leavevmode
 \epsfysize= 4cm 
\centerline{\epsfbox{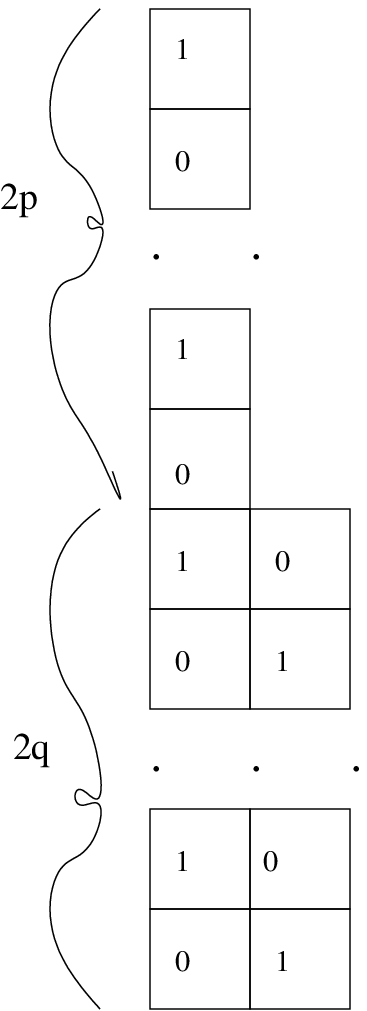}} 
\caption{\label{deca}}
\end{figure}

On a donc   
$n_0(D_{res}^{(1)}) =n_1(D_{res}^{(1)})=m \in {\N}$ et 
$n_0(D_{res}^{(2)}) =n_1(D_{res}^{(2)})=0$ i.e. l'ensemble
 $\{D_{res} / D\in {\bf  D_{\lambda,\mu}}\}$ est totalement ordonné.
Soit  $D,D' \in {\bf D_{\lambda,\mu}}$ on a $n_0(D^{(k)})
=n_0(D_{ind}^{(k)})+ n_0(D_{res}^{(k)})= n_0(D_{ind}^{'(k)})
+n_0(D_{res}^{(k)})$ et donc soit  $D\leq D'$ soit $D\geq D'$.

Supposons que $\lambda$ et $\mu$  comportent en commun une seule sous-partition de la forme $(k,k,k-1,k-1)$, $k\geq 1$.
On peut construire deux diagrammes gradués, l'un admettant comme sous-diagramme $D_1=(2k,2k,2k-2,2k-2)$, l'autre $D_2=(2k-1,2k-1,2k-1,2k-1)$. On ne s'intéresse qu'à ces sous-diagrammes les restes des deux diagrammes étant identiques.
Pour les comparer on efface les $2k-2$ premières colonnes qui sont identiques. Les deux diagrammes $D_1^{(2k-2)}$ et $D_2^{(2k-2)}$ restant sont ceux  de la figure 4 d'une part avec $p=2, q=0$ d'autre part avec $p=0, q=1$.  On en déduit que $D_{\scriptsize{\mbox{ max}(\lambda,\mu)}}$ existe.
\end{dem}

\begin{rem} - 
  Si $\lambda$ et $\mu$ comportent deux sous-partitions communes de la forme $(k,k,k-1,k-1)$ et $(l,l,l-1,l-1)$, $k\geq 1$ et $l\geq 1$, alors certains diagrammes ne sont pas comparables.
En effet, on construit les quatre sous-diagrammes possibles $D_1=(2k,2k,2k-2,2k-2, \ldots, 2l,2l,2l-2,2l-2)$,   $D_2=(2k,2k,2k-2,2k-2, \ldots, 2l-1,2l-1,2l-1,2l-1)$,  $D_3=(2k-1,2k-1,2k-1,2k-1, \ldots, 2l,2l,2l-2,2l-2)$ et  $D_4=(2k-1,2k-1,2k-1,2k-1, \ldots, 2l-1,2l-1,2l-1,2l-1)$.
On vérifie facilement que $D_2$ et $D_3$ ne sont pas comparables.
\end{rem}

\noi On en déduit le corollaire suivant:

\begin{coro}\label{coirr} Le cône nilpotent impair ${\mathcal N_1}$ de $\goth{ osp}(m,2n)$ est irréductible.
\end{coro}

\begin{dem} -  
 D'après la proposition \ref{irr}, il suffit  de donner tous les diagrammes gradués étiquetés possibles paramétrant une orbite maximale et de vérifier d'une part que dans aucun des cas $d_0$ n'est très paire d'autre part  que $\lambda$ et $\mu$  comportent au plus une seule sous-partition commune de la forme $(k,k,k-1,k-1)$, $k\geq 1$. Les diagrammes gradués étiquetés possibles sont décrits sur la figure \ref{orbmax}.

\begin{figure}[h!]
\leavevmode
 \epsfysize= 13cm 
{\epsfbox{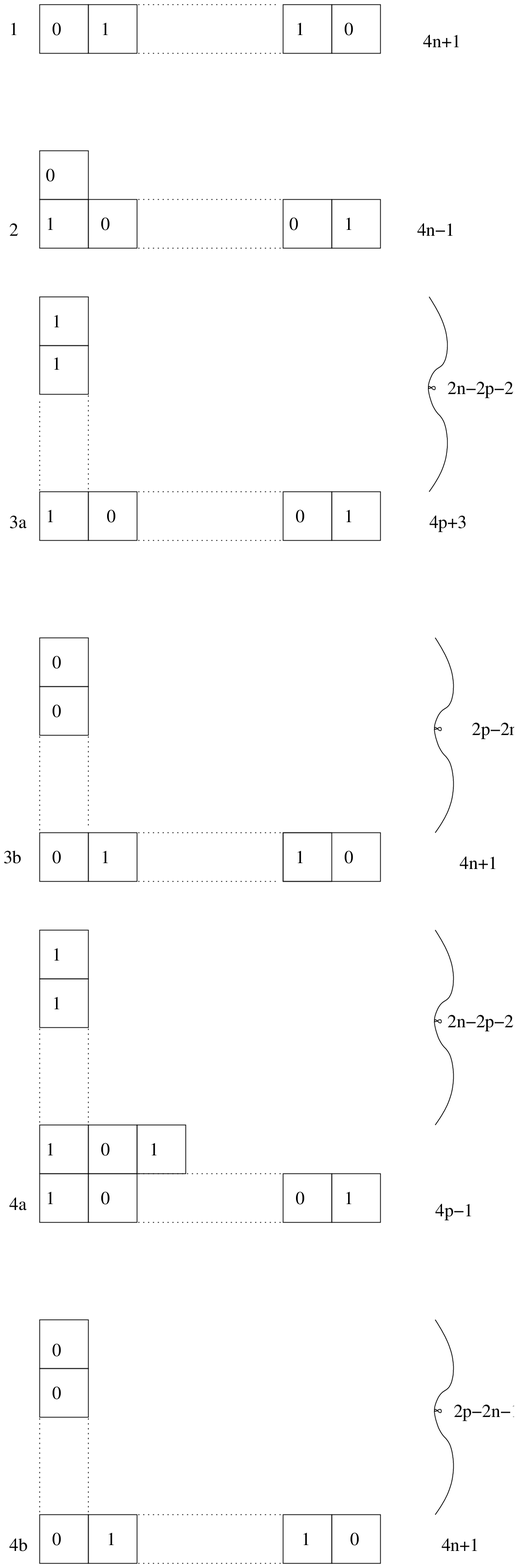}} 
\caption{\label{orbmax}}
\end{figure} 
 En effet, sauf dans le cas 4a, c'est un diagramme de forme "hook", dont la ligne la plus longue est impaire et de longueur maximale et la colonne est constitu\'ee de
cases sign\'ees isol\'ees en cardinal n\'ecessaire pour compl\'eter.
Dans le cas 4a  il est formé de deux lignes impaires et d'une colonne constitu\'ee de
cases sign\'ees isol\'ees en cardinal n\'ecessaire pour compl\'eter.

Il est évident que dans les cas ci-dessus  $\lambda$ et $\mu$ vérifient les hypothèses de la proposition 4.
\end{dem}

On s'intéresse ensuite à la fibre au dessus d'un élément quelconque de ${\mathcal N}_0$.
Pour en donner la dimension, on a besoin de la dimension des orbites nilpotentes impaires, que l'on 
trouve dans [{\bf KP}].

\begin{prop}{{\bf [ KP]}}
Soient $D\in {\mathcal D}(m,2n)$ un diagramme gradué et $O_D$ la $O(m)\times SP(2n)$-orbite associée.
On a 
$$\mbox {dim }O_D=\frac{1}{2}(\mbox{dim } \kappa(O_D)+\mbox{dim } V_0 \times \mbox{dim } V_1 -\Delta_D),$$
où \begin{displaymath}\Delta_D=\sum_{i\equiv 1\mod 2}P_i I_i,\end{displaymath} $P_i$ (resp. $I_i$) désignant le nombre de lignes 
paires (resp. impaires) de D de longueur $i$.
\end{prop}

On remarque que, dans le cas où $m$ est pair, la $O(m)\times SP(2n)$-orbite paramétrée par un diagramme 
$D$ tel que $d_0$ est une partition très paire, se décompose en deux composantes connexes de même dimension.

On peut à présent donner la dimension d'une fibre de $\kappa$ au dessus d'un élément :
\begin{prop}\label{pro3}  Soit $X_0\in {\mathcal N}_0$. Soit $(\lambda, \mu)\in 
P_{\goth{o}(m)etiq}\times P_{\goth{ sp}(2n)}$ tel que $X_0\in O_{\lambda, \mu}$ et $
D_{max(\lambda, \mu)}\in {\mathcal D}(m,2n)$ tel que $ O_{D_{max(\lambda, \mu)}}\subset \kappa^{-1}(O_{\lambda, \mu})$.
On a alors
$$\mbox{dim } \kappa^{-1}(X_0)=\frac{1}{2}(\mbox{dim } O_{\lambda, \mu}+\mbox{dim } V_0 \times \mbox{dim } V_1 -\Delta_{D_{max(\lambda, \mu)}}).$$
\end{prop}

\begin{dem} -  D'après la proposition \ref{irr} on sait qu'il existe une orbite maximale dans 
$\kappa^{-1}(O_{\lambda, \mu})$, notons-la $O_{D_{max(\lambda, \mu)}}$.
On regarde $\kappa: O_{D_{max(\lambda, \mu)}} \longrightarrow O_{\lambda, \mu}$.
Les fibres de ce morphisme sont équidimensionnelles et on a 
$\mbox{dim } \kappa^{-1}(X_0)= \mbox{dim } O_{\lambda, \mu}- \mbox{dim } O_{D_{max(\lambda, \mu)}}$.
Le résultat découle de la proposition précédente.
\end{dem}
On peut voir les dimensions des $G_0$-orbites nilpotentes impaires
de $\goth{ osp}(4,4)$ sur la  figure \ref{fig6}.

\section{Sous-alg\`ebres de Borel mixtes}

Rappelons d'abord les choix de bases que nous avons faits pour $V_0$ et $V_1$ :

$\bullet$ Si $dim(V_0) = 2p+1$, on choisit une base $e_0, \ldots , e_{2p}$ t.q. 
$B(e_i,e_j) = \delta _{i,2p-j}$, ce qui fait que les sous-espaces 
engendr\'es par $e_0, \ldots , e_{p-1}$ et $e_{p+1}, \ldots , e_{2p}$
sont totalement isotropes et en dualit\'e, le vecteur $e_p$ \'etant 
anisotrope.

$\bullet$ Si $dim(V_0) = 2p$, on choisit une base $e' _1, \ldots , e' _{2p}$ t.q. 
$B(e' _i,e' _j) = \delta _{i,2p+1-j}$, ce qui fait que les sous-espaces 
engendr\'es par $e' _1, \ldots , e' _{p}$ et $e' _{p+1}, \ldots , e' _{2p}$
sont totalement isotropes et en dualit\'e.

$\bullet$ On choisit une base $f_1, \ldots f_{2n}$ de $V_1$ telle que 
$B(f_i, f_j) = \delta _{i,2n+1-j}$ si $i\leq j$ et  
$B(f_i, f_j) = - \delta _{i,2n+1-j}$ si $i > j$. 
Les sous-espaces 
engendr\'es par $f_1, \ldots ,  f_{n}$ et $f _{n+1}, \ldots , f _{2n}$
sont totalement isotropes et en dualit\'e.

Soit $\goth b _0$ la sous-alg\`ebre de Borel de $\goth g_0$ 
qui pr\'eserve les drapeaux partiels d\'efinis par les suites de 
g\'en\'erateurs :

$\bullet$ Si $dim(V_0) = 2p+1$, $e_0, \ldots , e_{p-1}, f_1, \ldots ,  f_{n}$,

$\bullet$ Si $dim(V_0) = 2p$, $e' _1, \ldots , e' _{p}, f_1, \ldots ,  f_{n}$.

Cette sous-alg\`ebre pr\'eserve bien entendu les drapeaux complets 
d\'efinis par les bases choisies de $V_0$ et $V_1$, mais elle est d\'ej\`a
d\'etermin\'ee par les vecteurs sus-mentionn\'es. En effet, on compl\`ete
le drapeau en passant aux sous-espaces orthogonaux (on a 
$<e_0, \ldots , e_{p-1}> ^{\perp} = \-<e_0, \ldots , e_{p}>, 
<e' _1, \ldots , e' _{p}>^{\perp} = <e' _1, \ldots , e' _{p+1}>,
<f_1, \ldots ,  f_{n}>^{\perp} =<f_1, \ldots ,  f_{n}>,$ 

$<f_1, \ldots ,  f_{n-1}>^{\perp} 
=<f_1, \ldots ,  f_{n+1}> $, et ainsi de
suite).

On sait ({\bf{[Ka1], [DeL]}}) que les sous-alg\`ebres de Borel de $\goth g$ 
qui contiennent $\goth b _0$ ne sont pas deux-\`a-deux conjugu\'ees et que les classes de conjugaison 
sont en correspondance avec les ordres sur les vecteurs 
$e_0, \ldots , e_{p-1}, f_1, \ldots ,  f_{n}$ (resp. $e' _1, \ldots , 
e' _{p}, f_1, \ldots ,  f_{n}$) qui pr\'eservent l'ordre croissant sur
les $e_i$ (resp. $e' _i$) et les $f_j$.

Dans {\bf [Ka1]}, Kac choisit la sous-alg\`ebre de Borel correspondant
\`a l'ordre 

\noi $e_0, \ldots , e_{p-1}, f_1, \ldots ,  f_{n}$ (resp. $e' _1, \ldots , e' _{p}, f_1, \ldots ,  f_{n}$) pour \'etudier la th\'eorie des repr\'esentations. 

En observant la forme des matrices de repr\'esentants des orbites 
nilpotentes impaires, on constate que cette sous-alg\`ebre de Borel ne les rencontre pas toutes : elle ne rencontre que les orbites du c\^one autocommutant. 

Nous nous int\'eressons ici aux sous-alg\`ebres de Borel qui m\'elangent
le plus vecteurs pairs et impairs.

\begin{defi}\label{defmix}  La sous-alg\`ebre de Borel {\it mixte} contenant 
$\goth b _0$ est d\'efinie par l'ordre suivant sur les vecteurs :

\noi $\bullet$ Cas 1) $\goth{osp}(2n+1, 2n)$ : $e_0, f_1, e_1, \ldots , e_{n-1},
f_{n}$.

\noi $\bullet$ Cas 2) $\goth{osp}(2n, 2n)$ : $f_1, e' _1, f_2 \ldots ,
 e' _{n-1}, f_{n}, e' _n$.

\noi $\bullet$ Cas 3) $\goth{osp}(2p+1, 2n)$  avec

a) $p<n$ : 
$f_1, \ldots , f_{n-p}, e_0, f_{n-p+1}, e_1, \ldots , e_{p-1},
f_{n}$.

b) $p>n$ : 
$e_0, e_{p-n}, f_{1}, e_{p-n+1}, \ldots , e_{p-1},
f_{n}$.

\noi $\bullet$ Cas 4) $\goth{osp}(2p, 2n)$ avec

a) $p<n$ :
$f_1, \ldots , f_{n-p+1}, e' _1, f_{n-p+2}, e'_2, \ldots , f_n, e'_{p}$.

b) $p>n$ : 
$e'_1, e'_{p-n}, f_{1}, e'_{p-n+1}, \ldots , f_n, e'_{p}$.

Pour toute sous-alg\`ebre de Borel $\goth b_0$ de $\goth g _0$, on construit de la m\^eme mani\`ere une unique sous-alg\`ebre de Borel mixte de $\goth g$ contenant $\goth b_0$. On notera $\mathcal B _M$ l'ensemble des sous-alg\`ebres mixtes de $\goth g$.
\end{defi}

\begin{prop}  La sous-alg\`ebre de Borel mixte rencontre toutes les orbites nilpotentes impaires.
\end{prop}

\begin{dem} -  D'après le corollaire \ref{coirr} on sait qu'il n'existe qu'une orbite nilpotente maximale. Dans tous les cas le diagramme gradué étiqueté qui lui correspond est décrit sur la figure \ref{orbmax}.
On vérifie que cette orbite rencontre la sous-algèbre de Borel mixte.
On ne fait ici que le premier cas les autres se traitant de la même façon.
On le montre par récurrence sur $n$. Soit $(u, u^*)$
  un élément de  ${\mathcal N}_1$.  De $d_0$ (et donc des connaissances sur l'orbite nilpotente de $\goth{o}(2n+1)$) on déduit que le noyau de $u$ contient un vecteur isotrope de $V_0$, que l'on note $x$.  On considère l'hyperplan $H_0$ orthogonal à $x$ dans $V_0$.  On remarque que $u(H_0)\subsetneq V_1$ puisque $H_0$ contient $x$. Soit $y$ un élément de  $V_1$  orthogonal à
  $u(H_0)$. Soit $H_1$ l'espace vectoriel orthogonal à $y$ dans $V_1$.
  
On considère ensuite les espaces vectoriels
  
  $W_0 = H_0 / {\bf C}x$ and $W_1 = H_1 /{\bf C}y.$
  
  \noi On remarque que $W = W_0 \oplus W_1$  est un espace vectoriel orthosymplectique de 
  dimension $(2n-1,2n-2)$ et que $u$ induit un élément nilpotent impair de
  $\goth{osp}(W)$: on peut lui appliquer l'hypothèse de récurrence.
Et le résultat en découle.  
  
   Pour conclure, on aura besoin du lemme connu suivant :

\begin{leme}  Notons $\goth b _1$ la partie impaire d'une sous-alg\`ebre de Borel de $\goth g$. La r\'eunion des orbites sous $G_0$ qui rencontrent $\goth b_1$ est ferm\'ee dans $\goth g _1$.
\end{leme}

\begin{dem} - (du lemme)  Soit $Gr$ la grassmannienne des sous-espaces de dimension $mn$ de $\goth g _1$. Pour toute sous-alg\`ebre de Borel, $\goth b_1$ est un \'el\'ement de $Gr$ dont le stabilisateur dans $SO(m) \times Sp(2n)$ contient le sous-groupe $B_0$ de $G_0$ d'alg\`ebre de Lie la partie paire de la sous-alg\`ebre de Borel choisie. Son orbite dans $Gr$ est donc Zariski-ferm\'ee. Soit $I \subset Gr \times \goth g _1$ l'ensemble des couples $(V,u)$ dans $Gr \times \goth g _1$ avec $u \in V$.

On note $p_1, p_2$ les projections de $I$ sur chaque facteur. La r\'eunion des orbites qui rencontrent $\goth b _1$ est \begin{displaymath}\bigcup _{g \in SO(m) \times Sp(2n)} p_2(p_1 ^{-1}(g.\goth b _1)).\end{displaymath} Comme $p_1$ est continue et $p_2$ est propre, cet ensemble est bien Zariski-ferm\'e.
\end{dem}

La r\'eunion des orbites nilpotentes impaires qui rencontrent $\goth b_1$ est donc un ferm\'e qui contient l'orbite nilpotente impaire maximale, 
d'o\`u le r\'esultat.
\end{dem}

\section{D\'esingularisation du nilc\^one impair des super alg\`ebres de type
  $\goth{osp}(2n+1,2n)$ et $\goth{osp}(2n,2n)$.}
Nous commen\c cons par traiter en d\'etails le cas de $\goth{osp}(2n+1,2n)$.
 
On consid\`ere que, gr\^ace \`a la forme orthosymplectique, les
espaces $V_0$ et $V_1$ sont identifi\'es avec leurs duaux.

Soit $X = G_0 /B_0$ la vari\'et\'e des drapeaux de la partie paire de
$\goth g$. On note ${\mathcal O}_X$ le faisceau structural de $X$.

Notons
$$\{ 0 \} = E_0 \subset E_1 \subset \ldots \subset E_{2n+1} = V_0$$
le
drapeau tautologique de $V_0 \otimes {\mathcal O}_X$ et
$$\{ 0 \} = F_0 \subset F_1 \subset \ldots \subset F_{2n} = V_1$$
le
drapeau tautologique de $V_1 \otimes {\mathcal O}_X$.

Soit
$$\varphi _i : Hom (E_{i+1}, F_i) \longrightarrow Hom (E_{i}, F_i)$$
obtenue par restriction pour $1 \leq i \leq 2n$ et soit
$$\psi _i : Hom (E_{i}, F_{i-1}) \longrightarrow Hom (E_{i}, F_i)$$
obtenue par l'inclusion $F_{i-1} \subset F_i$, $2 \leq i \leq 2n$.

Soit $\chi : \oplus _{i} Hom (E_{i+1}, F_i) \longrightarrow \oplus
_{i} Hom (E_{i}, F_i)$ le morphisme dont la matrice par blocs est

$$\left ( \begin{array}{ccccc}
    \varphi _1 & 0 & 0 & \ldots & 0 \\
    -\psi _2 & \varphi _2 & 0 & \ldots & 0 \\
    0 &  -\psi _3 & \varphi _3 & \ldots & 0 \\
    \vdots & \vdots &\vdots &\vdots &\vdots \\
    0& 0 & 0 & \ldots & \varphi _{2n} \\
\end{array} \right ).$$
Alors $\chi$ est surjective donc $Ker \chi$ d\'efinit un fibr\'e au
dessus de $X$. On notera ${\tilde {\mathcal N}_1}$ ce fibr\'e.

\begin{prop} Le fibr\'e ${\tilde {\mathcal N}_1}$ au dessus de la vari\'et\'e de drapeaux $X$ est 
  une d\'e\-sin\-gu\-la\-ri\-sa\-tion du nilc\^one impair ${\mathcal N }_1$
  de $\goth g$
\end{prop}

\begin{dem} -  On a un morphisme $f : {\tilde {\mathcal N}_1} \longrightarrow \goth g _1$. Soit 
  ${\bf U}$ l'ouvert de ${\tilde {\mathcal N}_1}$ tel que : si une collection $(u_1, \ldots
  u_{2n})$ d'homomorphismes $u_i : E_{i+1} \longrightarrow F_i$ est un
  \'el\'ement de ${\tilde {\mathcal N}_1}$, $(u_1, \ldots u_{2n}) \in {\bf U}$ si et
  seulement si chaque $u_i$ est surjectif. On a $f(u_1, \ldots u_{2n})
  = u_{2n} \in \goth g_1$.
  
  Montrons que $u_{2n} \in {\mathcal N} _1$ : $u_1$ est la restriction de
  $u_{2n}$ \`a $E_2$. Sa restriction \`a $E_1$ est nulle. Donc le
  noyau de $u_{2n}$ contient $E_1$ et lui est donc \'egal puisque
  $u_{2n}$ est une application surjective. Le morphisme $u^* _{2n}$ a
  pour image $E_1 ^{\bot} = E_{2n}$ puique, pour tout \'el\'ement $g
  \in \goth g _1$, $(Ker g)^{\bot} = Im g^*$. Or $u_{2n}$ envoie $E_{2n}$ dans
  $F_{2n-1}$, on reproduit le raisonnement en rempla\c cant $V_0$ par
  $E_{2n} / E_1$ et $V_1$ par $F_{2n-1} / F_1$ : on se retrouve avec
  la m\^eme situation pour $\goth {osp}(2n-1,2n-2)$. On proc\`ede par
  r\'ecurrence : on v\'erifie que les choses marchent pour $\goth
  {osp} (3,2)$, o\`u on a la situation suivante:
  
  les drapeaux sont $$E_0 \subset E_1 \subset E_2 \subset E_3, F_0
  \subset F_1 \subset F_2$$
  et on a :
  
  $Ker u_2 = E_1$, $Im u_2 = F_2$, ${u_{2}} _{|E_2}: E_2 \rightarrow
  F_1$.
  
  D'autre part, $u_2 ^* : F_2 \rightarrow E_2 = E_1 ^{\bot}$.
  
  Donc $u_2 \circ u_2 ^*$ va de $F_2$ dans $F_1$. Or ${u_2 ^*}
  _{|F_1}$ va de $F_1$ dans $E_1$ ce qui fait que $u_2 ^* \circ u_2
  \circ u_2 ^*$ va de $F_2$ dans $E_1$, et $E_1$ est le noyau de $u_2$
  ce qui fait que $(u_2 \circ u_2 ^*)^2 =0$, donc $(u_2, u_2 ^*)$ est
  un \'el\'ement du nilc\^one ${\mathcal N}_1$.

  Le morphisme $f$ est un isomorphisme de ${\bf U}$ sur l'ouvert ${\bf
    S}$ de ${\mathcal N}_1$ form\'e des \'el\'ements surjectifs : si $(u,
  u^*) \in {\bf S}$, $dim \; Ker u =1$, donc $dim \; Ker (u^* \circ u)= 1$
   la suite des noyaux it\'er\'es des puissances de $u^* \circ u$ forme un drapeau
  complet et est donc un \'el\'ement de ${\bf U}$.
\end{dem}

\begin{rem} -  On peut r\'einterpr\'eter ${\tilde {\mathcal N}_1}$ comme suit :
$${\tilde {\mathcal N}_1} = \{(\goth b , x) \in \mathcal B _M \times 
\mathcal N _1 , x \in \goth b \}.$$
\end{rem}
\begin{theo}\label{th1}  Le fibr\'e vectoriel ${\tilde {\mathcal N}_1} = \{(\goth b , x) \in \mathcal B _M \times \mathcal N _1 , x \in \goth b \}$
au dessus de $G_0 /B_0$ est une d\'esingularisation de $\mathcal N _1$ pour
 $\goth{osp}(2n+1,2n)$ et $\goth{osp}(2n,2n)$.
\end{theo}

\begin{dem} -  Supposons que $\goth{ g}=\goth{osp}(2n+1,2n)$.
Soit $O_{D_{max}}$ la grosse orbite. Soit ${\bf U}=\{\mathcal B_M\times O_{D_{max}}\}$.
C'est un ouvert de $\widetilde{\mathcal N}_1$.

Soit $x=(u,v)\in O_{D_{max}}$. Montrons qu'on peut lui associer un 
drapeau complet $(W_i)$ correspondant à une sous-algèbre de Borel mixte.

On utilise le diagramme indexant l'orbite maximale de ${\mathcal N_1}$ 
(et donc la description des noyaux de $x\in O_{D_{max}}$ qu'il contient). 

On pose : 
\begin{itemize}
\item pour $i$ de $1$ à $n+1$, $W_{2i-1}=ker (x^i)_{\vert V_0}\oplus ker (x^{i-1})_{\vert V_1}$,  
\item pour  $j$ de $1$ à $n$, $W_{2j}=ker (x^{i})$.
\end{itemize}

D'après le diagramme de $O_{D_{max}}$  on remarque que
$\mbox{dim } W_{2i-1}= i+i-1$ de même $\mbox{dim } W_{2j}= \mbox{dim } ker(x^j)= 2j$.

De plus on a $W_i\cap V_0=E_{Ent(\frac{i+1}{2})}$, $W_i\cap V_1=F_{Ent(\frac{i}{2})}$, $((E_i),(F_j))\in X$.
On obtient donc que $(W_i)$ est un drapeau complet correspondant à une sous-algèbre de Borel mixte de $\goth{osp}(2n+1,2n)$.

Supposons à présent que $\goth{ g}=\goth{osp}(2n,2n)$.
Soit $O_{D_{max}}$ la grosse orbite. Soit ${\bf U}=\{\mathcal B_M\times O_{D_{max}}\}$.
C'est un ouvert de $\widetilde{\mathcal N}_1$.

Soit $x=(u,v)\in O_{D_{max}}$. Montrons qu'on peut lui associer un 
drapeau complet $(W_i)$.

On sait que pour $\goth{ o}(2n)$, $\lambda_{max}=(2n-1,1)$. On en déduit qu'il existe 
un vecteur $\alpha\in V_0$ anisotrope tel que $\alpha \in ker x$. Soit $V_0^{'}={V_0}_{\vert<\alpha>}$.
On pose :
\begin{itemize}
\item  pour $i$ de $1$ à $n$, $W_{2i-1}=ker (x^i)_{\vert V_1}\oplus ker (x^{i-1})_{\vert V_0^{'}}$, 
$W_{2i}=ker (x^{i})_{\vert (V_0^{'}\oplus V_1)}$,
\item pour $i$ de $n+1$ à $2n$, $W_{2i-1}=ker (x^{i-1})_{\vert V_0}\oplus ker (x^{i-1})_{\vert V_1}$, 
$W_{2i}=ker (x^{i})_{\vert V_1 }\oplus ker (x^{i-1})_{\vert V_0}$.
\end{itemize}

On a  $W_i\cap V_0=E_{Ent(\frac{i}{2})}$, $W_i\cap V_1=F_{Ent(\frac{i+1}{2})}$.
Comme ci-dessus on vérifie aisément que ce drapeau est complet. Et donc $(W_i)$ est bien un drapeau complet correspondant à une sous-algèbre de Borel mixte de $\goth{osp}(2n,2n)$.

\end{dem}

\begin{rem} -  On voit que cet \'enonc\'e est un analogue de la 
d\'esingularisation de Springer pour le c\^one nilpotent classique : si on note  $\widetilde{\mathcal N_0}=\{(\goth{ b}_0, x)\in X\times {\mathcal N_0}, \in { \goth b}_0)\}$, qui est le fibré cotangent à la variété des drapeaux  $X$, $\widetilde{\mathcal N_0}$ est une désingularisation du cône nilpotent pair ${\mathcal N_0}$.
\end{rem}
\begin{rem} -  Dans les cas o\`u les rangs des deux composantes de la partie paire diff\`erent, cette construction ne permet pas de reconstituer un drapeau complet de $V_0 \oplus V_1$ car tout \'el\'ement nilpotent de l'orbite maximale a un noyau de dimension au moins $2$.
\end{rem}

\begin{coro}

L'application  $g: \widetilde{\mathcal N_1}\longrightarrow \widetilde{\mathcal N_0}$ définie par 
$g((\goth{ b},(u,u^*)))=(\goth{ b}_0,(u^*\circ u,u\circ u^*))$
est un morphisme compatible avec  $\kappa$.
\end{coro}

\begin{dem} -  
Il suffit juste de montrer que le diagramme suivant commute :
$$\xymatrix{\\
 {\widetilde{\mathcal N_1}}\ar[r] \ar[d]^{\pi_1}&{\widetilde{\mathcal N_0}}\ar[d]^{\pi_0}&\\
{\mathcal N_1}\ar[r]^{\kappa}&{\mathcal N_0}&}.$$

Ce qui est évident car par construction on a
 $\pi_0(g((\goth{ b},(u,u^*)))=(u^*\circ u, u\circ u^*)$ et 
$\kappa((\pi_1(\goth{ b},(u, u^*)))=(u^*\circ u, u\circ u^*)$.
 \end{dem}

\section{Annexe}

Dans cette partie on rappelle quelques notions et notations (voir {\bf [Oh]}) sur les paires symétriques classiques, 
nécessaires à la démonstration de l'ordre sur les adhérences de Zariski des orbites 
nilpotentes impaires. 

Soit $G$ un groupe algébrique réductif  d'algèbre de Lie $\goth{ g}$ et $\theta$ une
 involution de $G$.
On note aussi $\theta$ l'involution induite sur $\goth{ g}$.
Soit $\goth{ g}=\goth{ k} + \goth{ p}$ la décompostion de Cartan de $\goth{ g}$ par 
rapport à $\theta$.
On note $(\goth{ g}, \goth{ k})$ la paire symétrique définie par $(G, \theta)$.
Soit $V$ un $\C$-espace vectoriel de dimension finie, $s$ une involution sur $V$. 
On munit $V$ d'une forme bilinéaire non dégénérée $\theta$-invariante. 
On pose $G(V):= \{g\in GL(V)/ g^*=g^{-1}\}$, $K(V):=\{g \in G(V)/ \theta(g)=g\}$, $\goth{ k}(V):= Lie K(V)$, 
$\goth{ p}(V)=\{X\in \goth{ g}(V)/ \theta(X)=-X\}$.
$(\goth{ g}(V), \goth{ k}(V))$ est une paire symétrique.

Les paires symétriques qui nous intéressent sont celles de type (DIII) c'est-à-dire 

$(\goth{ o}(2n, \C),\goth{ gl}(n,\C))$ et de type (CI) c'est-à-dire 
$(\goth{ sp}(2n,\C), \goth{ gl}(n, \C))$.
Les diagrammes indécomposables qui  paramètrent les $K(V)$-orbites nilpotentes dans $\goth{p}(V)$ sont les suivants.
Pour le type (DIII):
\begin{itemize}
\item deux lignes impaires de longueur $l$ paire,
\item deux lignes paires de longueur $l$ paire,
\item deux lignes, l'une paire l'autre impaire de longueur $l$ impaire.
\end{itemize}

Pour le type (CI):
\begin{itemize}
\item une ligne impaire de longueur paire,
\item une ligne paire de longueur paire,
\item deux lignes impaires de longueur $l$ impaire.
\end{itemize}


\section*{Bibliographie}

\noi {\bf [CM]} {\sc D. H. Collingwood, W. M. McGovern}, {\it
  Nilpotent orbits 
in semisimple Lie Algebras},
  Van Nostrand Reinhold, 1993.

\noi {\bf [deL]} {\sc J. van de Leur}, {\it Contragredient Lie superalgebras of
finite growth}, Utrecht University, 1986.

\noi {\bf [DL]} {\sc D. Z. Djokovi\'c, M. Litvinov}, {\it The closure ordering of nilpotent orbits 
of the complex symmetric pair $SO_{p+q}, SO_p\times SO_q$}, Canad. J. Math. {\bf 55},6 ,(2003), 1155-1190. 

\noi{\bf [DS]} {\sc M. Duflo, V. Serganova}, {\it On associated
  varieties for 
Lie Superalgebras },         
arxiv.math. RT/0507198.


\noi {\bf [Gru1]} {\sc C. Gruson}, {\it Finitude de l'homologie de
  certains 
modules
de dimension finie sur une super alg{\`e}bre de Lie}, Annales de
l'Institut Fourier, {\bf tome 47}, fascicule 2 (1997), pp. 531-553.

\noi {\bf [Gru2]} {\sc C. Gruson}, {\it Sur l'id{\'e}al autocommutant des
  super alg{\`e}bres de Lie basiques classiques et {\'e}tranges}, Annales de
l'Institut Fourier, {\bf tome 50}, fascicule 3 (2000), pp. 807-831.

\noi {\bf [Ja]} {\sc J. C. Jantzen}, {\it Nilpotent orbits in
  representation 
theory}, PIM, {\bf 228}, (2004), 1-211.  

\noi {\bf [Ka1]} {\sc V.G. Kac}, {\it Representations of classical Lie
superalgebras}, LNM 676, Springer (1978), pp. 597-626.
 
\noi {\bf [Ka2]} {\sc V.G. Kac}, {\it Some remarks on nilpotent orbits}, J.
of algebra, {\bf 64}, (1980), pp. 190-213.

\noi {\bf [KP]} {\sc H. Kraft, C. Procesi}, {\it On the geometry of conjugacy classes in classical groups}, Comment. Math. Helv. {\bf 57}, (1982), no. 4, pp. 539-602.

\noi {\bf [Oh]} {\sc T. Ohta}, {\it  The closures of nilpotent orbits
  in the 
classical symmetric pairs and their singularities }, Tohoku Math. J., 
{\bf 43}, (1991), pp. 161-211.


\noi {\bf [Sp1]} {\sc T. Springer}, {\it A construction of representations of
Weyl groups}, Invent. Math., {\bf  44}, (1978), pp. 279-293.

\noi {\bf [Sp2]} {\sc T. Springer}, {\it Quelques applications de la 
cohomologie
d'intersection}, S{\'e}minaire Bourbaki expos{\'e} {\bf 589}, F{\'e}vrier 1982.

\noi {\bf [Vu]} {\sc T . Vust}, {\it Sur la th{\'e}orie des invariants des groupes
classiques}, Ann.  Inst. Fourier, {\bf 26}, (1976), pp. 1-31.


\end{document}